\documentstyle[11pt,twoside]{article}
\title{\bf On the properties of the $(p,\nu)$-extension of the Whittaker function $M_{\kappa,\mu}(z)$}
\author{\sc S.A. Dar$^a$ and R.B. Paris$^b$\\
\\
${}^a\!$ {\em Department of Applied Sciences and Humanities, Faculty of Engineering }\\
{\em and Technology, Jamia Millia Islamia, New Delhi, 110025, India}\\ 
{\em E-Mail: showkatjmi34@gmail.com}\\
${}^b\!$ {\em Division of Computing and Mathematics, Abertay University,}\\
{\em Dundee DD1 1HG, UK}\\
{\em E-Mail: r.paris@abertay.ac.uk}
}
\begin{document}
\newcommand{\bee}{\begin{equation}}
\newcommand{\ee}{\end{equation}}
\newcommand{\br}{\biggr}
\newcommand{\bl}{\biggl}
\newcommand{\g}{\Gamma}
\def\f#1#2{\mbox{${\textstyle \frac{#1}{#2}}$}}
\def\dfrac#1#2{\displaystyle{\frac{#1}{#2}}}
\newcommand{\fr}{\frac{1}{2}}
\newcommand{\fs}{\f{1}{2}}
\date{}
\maketitle
\pagestyle{myheadings}
\markboth{\hfill \it S.A. Dar and R.B. Paris  \hfill}
{\hfill \it A $(p,\nu)$ extension of the Whittaker function \hfill}
\begin{abstract}
In this paper, we obtain a $(p,\nu)$-extension of the Whittaker function $M_{\kappa,\mu}(z)$ by using the extended confluent hypergeometric function of the first kind $\Phi_{p,\nu}(b;c;z)$ introduced in Parmar {\it et al.\/} [{\it J. Classical Anal.\/} {\bf 11} (2017) 81--106]. Also, we derive some of the main properties of this function, namely several integral representations, a summation formula, the analogue of Kummer's transformation formula, an asymptotic representation, the Mellin transform, a differential formula and some inequalities. 

\vspace{0.4cm}

\vspace{0.3cm}

\end{abstract}

\vspace{0.3cm}

\noindent $\,$\hrulefill $\,$

\vspace{0.2cm}

\begin{center}
{\bf 1. \  Introduction and Preliminaries}
\end{center}
\setcounter{section}{1}
\setcounter{equation}{0}
\renewcommand{\theequation}{\arabic{section}.\arabic{equation}}
In recent years several authors \cite{B,C,E1,E2,M,O1,PCP} have considered extensions of the familiar Gauss and confluent hypergeometric functions and carried out investigations into their basic properties. Since many special functions of interest in mathematical physics and in engineering applications, such as Bessel, Airy and Hermite functions, are particular cases of the confluent hypergeometric function, this last function occupies an important role in the theory of these special functions. It is hoped that generalisations of these hypergeometric functions may find useful application in various applied fields. 

The Gauss hypergeometric function $F(a,b;c;z)$ and the confluent hypergeometric function of the first kind $\Phi(b;c;z)$ (also denoted by $M(b,c,z)$) are given by
\begin{equation}\label{e4}
F(a,b;c;z)=\sum_{n=0}^\infty \frac{(a)_n(b)_n}{(c)_n}\,\frac{z^n}{n!}=\sum_{n=0}^{\infty}\frac{(a)_{n}B(b+n,c-b)}{B(b,c-b)}\frac{z^{n}}{n!}\quad(|z|<1),
\end{equation}
and
\begin{equation}\label{e5}
\Phi(b;c;z)=\sum_{n=0}^\infty\frac{(b)_n z^n}{(c)_n n!}=\sum_{n=0}^{\infty}\frac{B(b+n,c-b)}{B(b,c-b)}\frac{z^{n}}{n!}\qquad(|z|<\infty),
\end{equation}
where $(a)_n=\g(a+n)/\g(a)$ is the Pochhammer symbol, or rising factorial, and $B(x,y)$ is the Beta function defined by \cite[(5.12.1)]{DLMF}
\begin{equation}\label{e}
B(x,y)=\left\{\begin{array}{ll} \displaystyle{\int_{0}^{1}t^{x-1}(1-t)^{y-1}dt}, & (\Re(x)>0,  \Re(y)>0)\\
\\
\displaystyle{\frac{\Gamma(x)\Gamma(y)}{\Gamma(x+y)}}, & (x,y \neq 0, -1, -2, \ldots). \end{array}\right.
\end{equation}
The well-known integral representations of these functions are \cite[\S\S 13.4, 15.6]{DLMF}
\begin{equation}\label{e2}
F(a,b;c;z)=\frac{1}{B(b,c-b)}\int_{0}^{1}t^{b-1}(1-t)^{c-b-1}(1-zt)^{-a}dt,
\end{equation}
where $|\arg(1-z)|<\pi$, $\Re(c)>\Re(b)>0$, and
\begin{equation}\label{e3}
\Phi(b;c;z)=\frac{1}{B(b,c-b)}\int_{0}^{1}t^{b-1}(1-t)^{c-b-1} e^{zt}dt,
\end{equation}
where $\Re(c)>\Re(b)>0$.

The above-mentioned extensions of the hypergeometric functions have been achieved by employing suitable extensions of the Beta function. 
In 1997, Chaudhry {\it et al.} \cite[(1.7)]{B} introduced a $p$-extension of the Beta function given by
\begin{equation}\label{e7}
 B(x,y;p)=\int_{0}^{1}t^{x-1}(1-t)^{y-1} \exp\left(-\frac{p}{t(1-t)}\right) dt \qquad(\Re(p)>0).
\end{equation}
In a following paper \cite{C}, Chaudhry {\it et al.} employed this extended Beta function to define the extended Gauss hypergeometric function $F_{p}(a,b;c;z)$  and the extended confluent hypergeometric function $\Phi_{p}(b;c;z)$. These functions are defined by
\begin{equation}\label{e10}
F_{p}(a,b;c;z)
=\frac{1}{B(b,c-b)}\int_{0}^{1}t^{b-1}(1-t)^{c-b-1}(1-zt)^{-a}\exp\left(-\frac{p}{t(1-t)}\right)dt,
\end{equation}
where $|\arg(1-z)|<\pi$ and $\Re (p)>0$ ($p=0$, $\Re(c)>\Re(b)>0$), and
\begin{eqnarray}\label{e11}
\Phi_{p}(b;c;z)=\frac{1}{B(b,c-b)}\int_{0}^{1}t^{b-1}(1-t)^{c-b-1}\exp\left(zt-\frac{p}{t(1-t)}\right)dt,
\end{eqnarray}
where $\Re (p)>0$ ($p=0$, $\Re(c)>\Re(b)>0)$.

In \cite{PCP}, Parmar {\it et al.} further extended the Beta function $B(x,y;p)$ by adding one more parameter $\nu$ and considered the function
\begin{equation}\label{e13}
 B_{p,\nu}(x,y)=\sqrt{\frac{2p}{\pi}}\displaystyle{\int_{0}^{1}t^{x-\frac{3}{2}}(1-t)^{y-\frac{3}{2}}}
 K_{\nu+\frac{1}{2}}\left(\frac{p}{t(1-t)}\right) dt,
 \end{equation}
where $\Re(p)>0,\nu\geq0$ and $K_\nu(x)$ is the modified Bessel function of order $\nu$
These authors also introduced the extensions of $ F_{p}(a,b;c;z)$ and $\Phi_{p}(b;c;z)$ with the help of $B_{p,\nu}(x,y)$ given by:
\begin{equation}\label{e14}
F_{p,\nu}(a,b;c;z)\\
=\sum_{n=0}^{\infty}\frac{(a)_{n}B_{p,\nu}(b+n,c-b)}{B(b,c-b)}\frac{z^{n}}{n!}\qquad(|z|<1), 
\end{equation}
and
\begin{equation}\label{e15}
\Phi_{p,\nu}(b;c;z)=\sum_{n=0}^{\infty}\frac{B_{p,\nu}(b+n,c-b)}{B(b,c-b)}\frac{z^{n}}{n!}\qquad(|z|<\infty),
\end{equation}
where $\Re (p)>0$, $\nu\geq 0$.
Integral representations of the functions in (\ref{e14}) and (\ref{e15}) are given by \cite[\S 6]{PCP}
\begin{equation}\label{e16}
F_{p,\nu}(a,b;c;z)
=\frac{1}{B(b,c-b)}\sqrt{\frac{2p}{\pi}}\int_{0}^{1}t^{b-\frac{3}{2}}(1-t)^{c-b-\frac{3}{2}}
(1-zt)^{-a}K_{\nu+\frac{1}{2}}\left(\frac{p}{t(1-t)}\right)dt,
\end{equation}
where $|\arg(1-z)|<\pi$ and $\Re(p)>0$ ($p=\nu=0$, $\Re(c)>\Re(b)>0$), 
and
\begin{equation}\label{e17}
\Phi_{p,\nu}(b;c;z)=\frac{1}{B(b,c-b)}\sqrt{\frac{2p}{\pi}}\int_{0}^{1}t^{b-\frac{3}{2}}(1-t)^{c-b-\frac{3}{2}} e^{ zt}K_{\nu+\frac{1}{2}}\left(\frac{p}{t(1-t)}\right)dt,
\end{equation}
where $\Re(p)>0$ ($p=\nu=0$, $\Re(c)>\Re(b)>0$).

The Whittaker function of the first kind $M_{\kappa,\mu}(z)$  is defined by \cite[p.~334]{DLMF}, \cite[p.~337]{W}
\begin{equation}\label{e6}
M_{\kappa,\mu}(z)= z^{\mu+\frac{1}{2}}e^{-z/2}\Phi(\mu-\kappa+\fs;2\mu+1;z),
\end{equation}
provided $\mu$ ($2\mu\neq -1, -2, \ldots$),
where $\Phi(a;b;z)$ is the confluent hypergeometric function defined in (\ref{e5}).
An extension of the Whittaker function was given by Nagar {\it et al.} \cite[p.~74]{M} in the form
\begin{equation}\label{e12}
M_{p, k,\mu}(z)= z^{\mu+\frac{1}{2}}e^{-z/2}\Phi_{p}(\mu-k+\fs; 2\mu+1; z)\qquad (p\geq0)
\end{equation}
where $\Phi_p(\cdot)$ is defined in (\ref{e11}), which reduces to the ordinary Whittaker function when $p=0$.
A further extension of the Whittaker function has been made in \cite{DS}. Here the exponential factor appearing in $\Phi_p(z)$ is replaced by $\exp\,[-p/t-q/(1-t)]$, which reduces to the case in (\ref{e12}) when $p=q$.

In this paper, we consider the 
$(p,\nu)$-extension of the Whittaker function, which we denote by $M_{\kappa,\mu}^{(p,\nu)}(z)$, based on the extended confluent hypergeometric function $\Phi_{p,\nu}(z)$ defined in (\ref{e15}). Our work is  motivated in part by the following references (\cite{a2}, \cite{b1}, \cite{E1}, \cite{E2}, \cite{O1}). The $(p,\nu)$-extension of the Whittaker function is defined in Section 2 and some of its basic properties including integral representations and a differential property are given. In Section 3 we obtain the asymptotic behaviour of $M_{\kappa,\mu}^{(p,\nu)}(z)$
for $z\to+\infty$ and its Mellin transform, together with a related integral, in Section 4. In Section 5 we obtain some inequalities satisfied by our extended Whittaker function.
Some  concluding remarks are made in Section 6.

\vspace{0.6cm}

\begin{center}
{\bf 2. \  The extended Whittaker function and some basic properties}
\end{center}
\setcounter{section}{2}
\setcounter{equation}{0}
\renewcommand{\theequation}{\arabic{section}.\arabic{equation}}
The $(p,\nu)$-extension of the Whittaker function, which we denote by $M_{\kappa,\mu}^{(p,\nu)}(z)$, is defined 
in terms of the extended confluent hypergeometric function of the first kind $\Phi_{p,\nu}(z)$ in (\ref{e17}) by
\begin{equation}\label{e18}
M_{\kappa,\mu}^{(p,\nu)}(z)= z^{\mu+\frac{1}{2}} e^{-z/2}\Phi_{p, \nu}(\mu-k+\fs; 2\mu+1; z),
\end{equation}
where $\Re (p)\geq0$, $\nu\geq0$, $-\pi<\arg\,z\leq \pi$ and $\kappa$, $\mu$ ($2\mu\neq -1, -2 \ldots$) are complex parameters. If we set $p=\nu=0$, (\ref{e18}) reduces to the usual Whittaker function $M_{\kappa,\mu}(z)$.

The above extended Whittaker function satisfies the following differentiation formula:
\bee\label{e28}
\frac{d^n}{dz^n}\bl\{z^{-\mu-\frac{1}{2}} e^{z/2} M_{\kappa,\mu}^{(p,\nu)}(z)\br\}=\frac{(\mu-\kappa+\fs)_n}{(2\mu+1)_n}\, z^{-\mu-\frac{1}{2}-n/2}e^{z/2}\, M_{\kappa-\frac{1}{2}n,\mu+\frac{1}{2}n}^{(p,\nu)}(z)
\ee
for $n=0, 1, 2, \ldots\ $. This follows immediately from the result \cite[(6.9)]{PCP}
\[\frac{d^{n}}{dz^{n}}\left\{\Phi_{p,\nu}(b;c;z)\right\}=\frac{(b)_{n}}{(c)_{n}}\Phi_{p,\nu}(b+n;c+n;z)\]
for non-negative integer $n$,
upon noticing that the quantity on the left-hand side of (\ref{e28}) is, from (\ref{e18}), the $n$th derivative of the function $\Phi_{p,\nu}(\mu-\kappa+\fs;2\mu+1;z)$.

An integral representation follows from (\ref{e18}) and (\ref{e17}) in the form
\[M_{\kappa,\mu}^{(p,\nu)}(z)= \frac{z^{\mu+\frac{1}{2}} e^{-z/2}}{B(\mu-\kappa+\fs, \mu+\kappa+\fs)}\hspace{4cm}\] 
\[\hspace{3cm}\times
\sqrt{\frac{2p}{\pi}}\int_{0}^{1}t^{\mu-\kappa-1}(1-t)^{\mu+\kappa-1} e^{zt} K_{\nu+\frac{1}{2}}\left(\frac{p}{t(1-t)}\right)dt,\]
\bee\label{e20}
(\nu\geq0,\quad \Re (p)>0;\ \  p=0 \ \ \mbox{and}\ \ \Re (\mu\pm\kappa+\fs)>0).
\ee
 Alternative representations can be obtained by making the change of variable $t=(u-\alpha)/(\beta-\alpha)$, with $\beta-\alpha>0$, to yield
\[M_{\kappa,\mu}^{(p,\nu)}(z)=                                                                                                           \frac{(\beta-\alpha)^{-2\mu+1}z^{\mu+\frac{1}{2}} e^{-z/2}}{B(\mu-\kappa+\frac{1}{2},                                                                     \mu+\kappa+\frac{1}{2})}
\sqrt{\frac{2p}{\pi}}\int_{\alpha}^{\beta}(u-\alpha)^{\mu-\kappa-1}(\beta-u)^{\mu+\kappa-1}\]
\bee\label{e21} \times \exp\left\{\frac{z(u-\alpha)}{\beta-\alpha}\right\}K_{\nu+\frac{1}{2}}\left\{\frac{p(\beta-\alpha)^{2}}{(u-\alpha)(\beta-u)}\right\}du,
\ee
and $t=u/(1+u)$ to yield
\[M_{\kappa,\mu}^{(p,\nu)}(z)=                                                                                                           \frac{z^{\mu+\frac{1}{2}} e^{-z/2}}{B(\mu-\kappa+\frac{1}{2},                                                                     \mu+\kappa+\frac{1}{2})} \hspace{6cm}\]
\bee\label{e23}
\hspace{3cm}\times \sqrt{\frac{2p}{\pi}}\int_{0}^{\infty}\frac{u^{\mu-\kappa-1}}{(1+u)^{2\mu}} \exp\left\{\frac{zu}{1+u}\right\}K_{\nu+\frac{1}{2}}\left\{\frac{p(1+u)^{2}}{u}\right\}du,
\ee
valid under the same conditions as (\ref{e20}).

A special case of (\ref{e21}) is obtained by letting $\alpha=-1$, $\beta=1$ to find
\[M_{\kappa,\mu}^{(p,\nu)}(z)=                                                                                                           \frac{2^{-2\mu+1}z^{\mu+\frac{1}{2}}}{B(\mu-\kappa+\frac{1}{2},                                                                     \mu+\kappa+\frac{1}{2})}\hspace{6cm}\]
\[
\times \sqrt{\frac{2p}{\pi}}
\int_{-1}^{1}(1+u)^{\mu-\kappa-1}(1-u)^{\mu+\kappa-1} e^{zu/2}K_{\nu+\frac{1}{2}}\left(\frac{4p}{1-u^2}\right)du,\]
\bee\label{e22}
(\nu\geq0,\quad \Re (p)>0;\ \  p=0 \ \ \mbox{and}\ \ \Re (\mu\pm\kappa+\fs)>0).
\ee
When $\nu=0$, we have upon making use of the result $K_\frac{1}{2}(z)=(\pi/2z)^{1/2}e^{-z}$, the representation
\[M_{\kappa,\mu}^{(p,0)}(z)=                                                                                                           \frac{2^{-2\mu}z^{\mu+\frac{1}{2}}}{B(\mu-\kappa+\frac{1}{2},                                                                     \mu+\kappa+\frac{1}{2})}\hspace{5cm}\]
\[\times 
\int_{-1}^{1}(1+u)^{\mu-\kappa-\frac{1}{2}}(1-u)^{\mu+\kappa-\frac{1}{2}}\exp\left\{\frac{zu}{2}-\frac{4p}{(1-u^2)}\right\}du,\]
\bee\label{E23}
(\nu\geq0,\quad \Re (p)>0;\ \  p=0 \ \ \mbox{and}\ \ \Re (\mu\pm\kappa+\fs)>0).
\ee
When $p=0$, this last result  reduces to the well-known representation for the usual Whittaker function $M_{\kappa,\mu}(z)$ in \cite[(13.16.1)]{DLMF} valid for $\Re (\mu\pm\kappa+\fs)>0$.

\newtheorem{theorem}{Theorem}
\begin{theorem}$\!\!\!.$\ When $\nu=n$, where $n$ is a non-negative integer, the following summation formula holds:
\[M_{\kappa,\mu}^{(p,n)}(z)=z^{\mu+\frac{1}{2}}e^{-z/2} \sum_{k=0}^n\frac{(n+k)! (2p)^{-k}}{(n-k)! k!}\,\frac{(\mu-\kappa+\fs)_k (\mu+\kappa+\fs)_k}{(2\mu+1)_{2k}}\hspace{4cm}\]
\bee\label{e51}
\hspace{4cm}\times\Phi_p(\mu-\kappa+k+\fs;2\mu+2+1;z),
\ee
where $\Re (p)>0$.
\end{theorem}
\noindent{\bf Proof.}\ \ When $\nu=n$, $n=0, 1, 2, \ldots\,$, the modified Bessel function $K_{n+\frac{1}{2}}(z)$ has the expansion \cite[(10.49.12)]{DLMF}
\[K_{n+\frac{1}{2}}(z)=\sqrt{\frac{\pi}{2z}}\, e^{-z} \sum_{k=0}^n \frac{a_k(n)}{(2z)^k},\qquad a_k(n):=\frac{(n+k)!}{(n-k)! k!}.\]
It then follows from (\ref{e20}) that
\begin{eqnarray*}
M_{\kappa\mu}^{(p,n)}(z)&=&\frac{z^{\mu+\frac{1}{2}}e^{-z/2}}{B(\mu-\kappa+\fs,\mu+\kappa+\fs)}\sum_{k=0}^n\frac{a_k(n)}{(2p)^k}\\
&&\hspace{3cm}\times\int_0^1t^{\mu-\kappa+k-\frac{1}{2}}(1-t)^{\mu+\kappa+k-\frac{1}{2}} \,\exp\,\bl(zt-\frac{p}{t(1-t)}\br) dt\\
&=&z^{\mu+\frac{1}{2}}e^{-z/2} \sum_{k=0}^n\frac{a_k(n)}{(2p)^k}\,\frac{B(\mu-\kappa+k+\fs,\mu+\kappa+k+\fs)}{B(\mu-\kappa+\fs,\mu+\kappa+\fs)}\\
&&\hspace{3cm}\times\Phi_p(\mu-\kappa+k+\fs;2\mu+2k+1;z)
\end{eqnarray*}
upon evaluation of the integral by means of (\ref{e11}). Some straightforward  manipulation then yields the result stated in (\ref{e51}).

To conclude this section we establish the analogue of the Kummer transformation for $M_{\kappa,\mu}^{(p,\nu)}(z)$.
\begin{theorem}$\!\!\!.$\ The following Kummer-type transformation holds:
\bee\label{e25}
z^{-\mu-\frac{1}{2}} M_{\kappa,\mu}^{(p,\nu)}(z)=(-z)^{\mu-\frac{1}{2}} M_{-\kappa,\mu}^{(p,\nu)}(-z),
\ee
for $p\geq 0$, $\nu\geq 0$ and $2\mu\neq -1, -2, \ldots\ $.
\end{theorem}
\noindent {Proof.}\ \ From the result \cite[(6.11)]{PCP}
\[\Phi_{p,\nu}(b;c;z)= e^z\,\Phi_{p,\nu}(c-b;c;-z),\]
we have, provided $2\mu$ is not a negative integer,
\begin{eqnarray*}
M_{\kappa\mu}^{(p,\nu)}(-z)&=&(-z)^{\mu+\frac{1}{2}}\, e^{z/2}\, \Phi_{p,\nu}(\mu-\kappa+\fs;2\mu+1;-z)\\
&=&(-z)^{\mu+\frac{1}{2}}\, e^{-z/2}\, \Phi_{p,\nu}(\mu+\kappa+\fs;2\mu+1;z).
\end{eqnarray*}
Application of (\ref{e18}) again to the right-hand side of this last expression then yields the result stated in (\ref{e25}).

When $p=\nu=0$, (\ref{e25}) reduces to the standard transformation formula for $M_{\kappa,\mu}(z)$ given in \cite[p.~338]{W}.
\vspace{0.6cm}

\begin{center}
{\bf 3. \  The asymptotic behaviour of $M_{\kappa,\mu}^{(p,\nu)}(z)$ for $x\to+\infty$}
\end{center}
\setcounter{section}{3}
\setcounter{equation}{0}
\renewcommand{\theequation}{\arabic{section}.\arabic{equation}}
We suppose that $p>0$, $\nu\geq0$ and consider positive values of $z$ ($=x$). Then the following theorem holds:
\begin{theorem}$\!\!\!.$\ For $p>0$, $\nu\geq0$, the asymptotic behaviour of the extended Whittaker function is independent of $\nu$ to leading order and is given by
\bee\label{e330}
M_{\kappa,\mu}^{(p,\nu)}(x)\sim A x^{(\mu-\kappa)/2} \,\exp\,(\fs x-2\sqrt{px})\qquad (x\to+\infty),
\ee
where 
\[A=\frac{\sqrt{\pi} e^{-p} p^{(\mu+\kappa)/2}}{B(\mu-\kappa+\fs,\mu+\kappa+\fs)}.\]
\end{theorem}
\noindent{\bf Proof.}\ \ If we make the substitution $t\to 1-t$ in the integral in (\ref{e20}) we find
\[M_{\kappa,\mu}^{(p,\nu)}(x)=\frac{x^{\mu+\frac{1}{2}}e^{x/2}}{B(\mu-\kappa+\fs,\mu+\kappa+\fs)}\sqrt{\frac{2p}{\pi}}\hspace{6cm}\]
\[\hspace{4cm}\int_0^1t^{\mu+\kappa-1}(1-t)^{\mu-\kappa-1} e^{-xt}\,K_{\nu+\frac{1}{2}}\bl(\frac{p}{t(1-t)}\br)dt.\]
As $x\to+\infty$, the main contribution to the integral comes from the neighbourhood of $t=0$, where the Bessel function may be replaced by its asymptotic form \cite[(10.25.3)]{DLMF}
\bee\label{e333}
K_{\nu+\frac{1}{2}}\bl(\frac{p}{t(1-t)}\br)\sim \sqrt{\frac{\pi}{2p}}\,\{t(1-t)\}^{1/2}\,\exp\,\bl(-\frac{p}{t}-\frac{p}{1-t}\br)\qquad (t\to 0).
\ee
This yields the approximate representation
\bee\label{e331}
M_{\kappa,\mu}^{(p,\nu)}(x)\simeq\frac{x^{\mu+\frac{1}{2}}e^{x/2}}{B(\mu-\kappa+\fs,\mu+\kappa+\fs)}\,I(x),
\ee
where
\[I(x)=\int_0^1 \exp\,\bl(-xt-\frac{p}{t}\br) t^{\mu+\kappa-\frac{1}{2}}(1-t)^{\mu-\kappa-\frac{1}{2}} e^{-p/(1-t)}dt.\]

With the rescaling of the variable $t=\tau/\sqrt{x}$, we obtain
\[I(x)=x^{-\frac{1}{2}(\mu+\kappa+\frac{1}{2})} \int_0^{\sqrt{x}} e^{-\sqrt{x} \psi(\tau)} f(\tau)\,d\tau,\]
where
\[\psi(\tau)=\tau+p/\tau,\qquad f(\tau)=\tau^{\mu+\kappa-\frac{1}{2}} e^{-p} \{1+O(x^{-1/2})\}.\]
The phase function $\psi(\tau)$ has a saddle point where $\psi'(\tau)=0$; that is, at the point $\tau_s=\sqrt{p}$. It is easily verified that the integration path in $I(x)$ is the steepest descent path through the saddle $\tau_s$. Application of the saddle-point method \cite[(2.4.15)]{DLMF} as $x\to+\infty$ then shows that 
\begin{eqnarray*}
I(x)&=&\frac{2e^{-\sqrt{x}\psi(\tau_s)}}{x^{\frac{1}{2}(\mu+\kappa+\frac{1}{2})}}\,\sqrt{\frac{\pi}{2\sqrt{x} \psi''(\tau_s)}}\,f(\tau_s) \{1+O(x^{-1/2})\}\\
&=&\frac{e^{-2\sqrt{px}}}{x^{\frac{1}{2}(\mu+\kappa+1)}}\,\sqrt{\pi} e^{-p} p^{\frac{1}{2}(\mu+\kappa)}
\{1+O(x^{-1/2})\},
\end{eqnarray*}
where $\psi(\tau_s)=2\sqrt{p}$ and $\psi''(\tau_s)=2/\sqrt{p}$.
Insertion of the above estimate for $I(x)$ in (\ref{e331}) then leads to the result stated in (\ref{e330}).

It is possible to extend the above analysis to show that the leading asymptotic form in (\ref{e330}) holds for $|\arg\,z|<\fs\pi$; we omit these details.

\vspace{0.6cm}

\begin{center}
{\bf 4. \  Mellin transforms of $M_{\kappa,\mu}^{(p,\nu)}(z)$}
\end{center}
\setcounter{section}{4}
\setcounter{equation}{0}
\renewcommand{\theequation}{\arabic{section}.\arabic{equation}}
The Mellin transform of a locally integrable function $f(x)$ on $(0,\infty)$ is given by (see, for example,  \cite[p.29]{DLMF})
\begin{equation}\label{e30}
{\cal M}\left\{f(x)\right\}(s)=\int_{0}^{\infty}x^{s-1}f(x)\,dx
\end{equation}
which defines an analytic function in its strip of analyticity. 

\begin{theorem}$\!\!\!.$\ The following Mellin transform holds:
\[\int_0^\infty p^{s-1} M_{\kappa,\mu}^{(p,\nu)}(z)\,dp=\frac{2^{s-1} z^{\mu+\frac{1}{2}} e^{-z/2}}{\sqrt{\pi}\,B(\mu-\kappa+\fs,\mu+\kappa+\fs)}\,\Gamma\bl(\frac{s-\nu}{2}\br) \Gamma\bl(\frac{s+\nu+1}{2}\br)\]
\bee\label{e31}
\times B(\mu-\kappa+s+\fs, \mu+\kappa+s+\fs)\,\Phi(\mu-\kappa+s+\fs; 2\mu+2s+1; z),
\ee
where $\Re (s)>\nu$ and $\Re (\mu\pm\kappa+\fs)>0$.
\end{theorem}
\noindent{\bf Proof.}\ \ From the fact that $K_{\nu+\frac{1}{2}}(t)=O(t^{-\nu-\frac{1}{2}})$ as $t\to 0$ and from the
asymptotic form (\ref{e333}), it is easily seen from (\ref{e18}) and the integral (\ref{e20}) that
\[M_{\kappa,\mu}^{(p,\nu)}(z)=\left\{\begin{array}{ll}O(p^{-\nu}) & p\to 0\\
\\
O(e^{-4p}) & p\to+\infty.\end{array}\right.\]
The strip of analyticity of the integral on the left-hand side of (\ref{e31}) is therefore $\Re (s)>\nu$.

From (\ref{e18}) and the definition of $\Phi_{p,\nu}(z)$ in (\ref{e17}) we obtain
\[
\int_0^\infty p^{s-1} M_{\kappa,\mu}^{(p,\nu)}(z)\,dp=z^{\mu+\frac{1}{2}} e^{-z/2} \int_0^\infty p^{s-1} \Phi_{p,\nu}(\mu-\kappa+\fs;2\mu+1;z)\,dp\]
\[=\frac{z^{\mu+\frac{1}{2}} e^{-z/2}}{B(\mu-\kappa+\fs, \mu+\kappa+\fs)}\,\sqrt{\frac{2}{\pi}}\hspace{6cm} \]
\[\hspace{3cm}\times\int_0^1 t^{\mu-\kappa-1}(1-t)^{\mu+\kappa-1}\,e^{zt}\bl\{\int_0^\infty p^{s-\frac{1}{2}} K_{\nu+\frac{1}{2}}\bl(\frac{p}{t(1-t)}\br)\br\}dt\]
upon reversal of the order of integration. Application of the result \cite[(10.43.19)]{DLMF}
\[\int_0^\infty w^{s-\frac{1}{2}}K_{\alpha+\frac{1}{2}}(w)\,dw=2^{s-\frac{3}{2}} \Gamma\bl(\frac{s-\alpha}{2}\br) \Gamma\bl(\frac{s+\alpha+1}{2}\br)\qquad (|\Re (\alpha)|<\Re (s))\]
after the change of variable $w=p/(t(1-t))$ in the inner integral produces
\[\int_0^\infty p^{s-1} M_{\kappa,\mu}^{(p,\nu)}(z)\,dp=\frac{2^{s-1} z^{\mu+\frac{1}{2}} e^{-z/2}}{\sqrt{\pi}\,B(\mu-\kappa+\fs,\mu+\kappa+\fs)}\,\Gamma\bl(\frac{s-\nu}{2}\br) \Gamma\bl(\frac{s+\nu+1}{2}\br)\]
\[\times\int_0^1 t^{\mu-\kappa+s-\frac{1}{2}}(1-t)^{\mu+\kappa+s-\frac{1}{2}} e^{zt}dt.\]
Identification of the above integral as a confluent hypergeometric function by (\ref{e3}) then yields the right-hand side of (\ref{e31}) when $\Re (s)>\nu$ and $\Re (\mu\pm\kappa+\fs)>0$.

\begin{theorem}$\!\!\!.$\ Let $p>0$, $\nu\geq 0$ and $\alpha>0$, $\beta>0$ be arbitrary parameters satisfying $2\alpha-\beta\geq 0$. Define the variable $\chi:=2\beta/(2\alpha+\beta)$, so that $0<\chi\leq 1$. Then
\[\int_0^\infty x^{\alpha-1}e^{-\alpha x}  M_{\kappa,\mu}^{(p,\nu)}(\beta x)\,dx\hspace{8cm}\]
\bee\label{e32}
\hspace{1cm}=\beta^{-\alpha} \Gamma(\mu+\alpha+\fs)\,\chi^{\mu+\alpha+\frac{1}{2}}\,F_{p,\nu}(\mu+\alpha+\fs, \mu-\kappa+\fs;2\mu+1;\chi),
\ee
provided $\Re (\mu+\alpha+\fs)>0$.
\end{theorem}
\noindent{\bf Proof.}\ \ From (\ref{e20}), it is evident that $M_{\kappa,\mu}^{(p,\nu)}(x)=O(x^{\mu+\frac{1}{2}})$ as $x\to 0$ and, from (\ref{e330}), we have
\[M_{\kappa,\mu}^{(p,\nu)}(\beta x)=O(x^{(\mu-\kappa)/2} \exp\,[\fs\beta x-2\sqrt{p\beta x}])\qquad (x\to+\infty).\]
Hence the integral on the left-hand side of (\ref{e32}) converges for $\Re (\mu+\alpha+\fs)>0$ and $2\alpha-\beta\geq0$. 

From (\ref{e18}) and (\ref{e20}) we have
\[I\equiv\int_0^\infty \!\!x^{\alpha-1}e^{-\alpha x}M_{\kappa,\mu}^{(p,\nu)}(\beta x)\,dx
=\beta^{\mu+\frac{1}{2}}\int_0^\infty\!\! x^{\alpha+\mu-\frac{1}{2}} e^{-(2\alpha+\beta)x/2}\,\Phi_{p,\nu}(\mu-\kappa+\fs;2\mu+1;\beta x)\,dx\]
\[=\frac{\beta^{\mu+\frac{1}{2}}}{B(\mu-\kappa+\fs,\mu+\kappa+\fs)}\,\sqrt{\frac{2p}{\pi}}\]
\[\times\int_0^1 t^{\mu-\kappa-1}(1-t)^{\mu+\kappa-1} K_{\nu+\frac{1}{2}}\bl(\frac{p}{t(1-t)}\br)
 \bl\{\int_0^\infty x^{\alpha+\mu-\frac{1}{2}} \exp\,[-(2\alpha+\beta)x(1-\chi t)]\,dx\br\}dt,\]
where $\chi=2\beta/(2\alpha+\beta)$ and we have interchanged the order of integration. Provided $0<\chi<1$ and $\Re (\mu+\alpha+\fs)>0$, the inner integral can be evaluated by Euler's integral formula $\int_0^\infty x^{\rho-1}e^{-\gamma x}dx=\gamma^{-\rho}\g(\rho)$, $\gamma>0$, to find
\[I=\frac{\chi^{\mu+\alpha+\frac{1}{2}} \g(\mu+\alpha+\fs)}{\beta^\alpha B(\mu-\kappa+\fs,\mu+\kappa+\fs)}\,\sqrt{\frac{2p}{\pi}}\hspace{6cm}\]
\[\hspace{4cm}\times \int_0^1 t^{\mu-\kappa-1}(1-t)^{\mu+\kappa-1}(1-\chi t)^{-\mu-\alpha-\frac{1}{2}} K_{\nu+\frac{1}{2}}\bl(\frac{p}{t(1-t)}\br)\,dt.\]
Finally, identifying the above integral as the extended hypergeometric function $F_{p,\nu}(z)$ in (\ref{e16}) we obtain the right-hand side of (\ref{e32}).

The result (\ref{e32}) has been established for $0<\chi<1$ ($2\alpha-\beta>0$). As both sides of this equation are analytic functions when $\chi=1$, the result (\ref{e32}) holds by analytic continuation for $0<\chi\leq 1$; that is, for $2\alpha-\beta\geq 0$.
\bigskip

\noindent{\bf Corollary.}\ \ If we set $\beta=2\alpha$ in (\ref{e32}), we find
\[\int_0^\infty x^{\alpha-1}e^{-\alpha x} M_{\kappa,\mu}^{(p,\nu)}(2\alpha x)\,dx=\frac{\g(\mu+\alpha+\fs)}{(2\alpha)^\alpha}\,F_{p,\nu}(\mu+\alpha+\fs, \mu-\kappa+\fs;2\mu+1;1).\]
Use of the evaluation \cite[Eq.~(6.13)]{PCP}
\[F_{p,\nu}(a,b;c;1)=\frac{B_{p,\nu}(b,c-a-b)}{B(b,c-b)}\qquad (p>0),\]
then shows that
\bee\label{e33}
\int_0^\infty x^{\alpha-1}e^{-\alpha x} M_{\kappa,\mu}^{(p,\nu)}(2\alpha x)\,dx=\frac{\g(\mu+\alpha+\fs)}{(2\alpha)^\alpha}\,\frac{B_{p,\nu}(\mu-\kappa+\fs, \kappa-\alpha)}{B(\mu\!-\!\kappa\!+\!\fs,\mu\!+\!\kappa\!+\!\fs)}
\ee
for $p>0$, $\nu\geq0$, $\alpha>0$ and $\Re (\mu+\alpha+\fs)>0$.
\vspace{0.6cm}

\begin{center}
{\bf 5. \  Inequalities for $M_{\kappa,\mu}^{(p,\nu)}(z)$}
\end{center}
\setcounter{section}{5}
\setcounter{equation}{0}
\renewcommand{\theequation}{\arabic{section}.\arabic{equation}}
\begin{theorem}$\!\!\!.$\ Let $\Re (p)>0$, $\nu>0$ and the parameters $\kappa$, $\mu$ be real and such that $\mu\pm\kappa+\fs>0$. Then the following upper bound holds:
\bee\label{e55}
|M_{\kappa,\mu}^{(p,\nu)}(z)|<C |z|^{\mu+\frac{1}{2}} e^{-\Re (z)/2}\,\frac{B(\mu\!-\!\kappa\!+\!\nu\!+\!\fs,\mu\!+\!\kappa\!+\!\nu\!+\!\fs)}{B(\mu\!-\!\kappa\!+\!\fs,\mu\!+\!\kappa\!+\!\fs)}\,\Phi(\mu-\kappa+\nu+\fs;2\mu+2\nu+1;\Re (z)),
\ee
where 
\[C=\frac{2^\nu |p|^{\nu+\frac{1}{2}} \g(\nu+\fs)}{\sqrt{\pi}\,(\Re (p))^{2\nu+1}}.\]
\end{theorem}
\noindent{\bf Proof.}\ \ From \cite[(5.3)]{DP1}, we have the bound for the modified Bessel function given by
\[|K_{\nu+\frac{1}{2}}(\zeta)|<\frac{1}{2}\bl(\frac{2|\zeta|}{(\Re (\zeta))^2}\br)^{\nu+\frac{1}{2}} \g(\nu+\fs)\qquad(\nu>0,\ \Re (\zeta)>0).\]
With $\zeta=p/(t(1-t))$ and $t\in (0,1)$ we therefore find, for $\Re (p)>0$,
\[
\left| K_{\nu+\frac{1}{2}}\left(\frac{p}{t(1-t)}\right)\right|<\frac{1}{2}\left(\frac{2|p|t(1-t)}{(\Re (p))^{2}}\right)^{\nu+\frac{1}{2}}\Gamma(\nu+\fs).
\]
For ease of presentation we shall suppose that $\kappa$ and $\mu$ are real and that $\mu\pm\kappa+\fs>0$. Then, from (\ref{e20}), it follows that
\begin{eqnarray*}
&&|M_{\kappa\mu}^{(p,\nu)}(z)|\\
\!\!&<&\!\!\frac{|z|^{\mu+\frac{1}{2}}e^{-\Re (z)/2}}{B(\mu\!-\!\kappa\!+\!\fs,\mu\!+\!\kappa\!+\!\fs)}\sqrt{\frac{2|p|}{\pi}}\int_0^1 t^{\mu-\kappa-1}(1-t)^{\mu+\kappa-1}e^{\Re (z)t}\left| K_{\nu+\frac{1}{2}}\left(\frac{p}{t(1-t)}\right)\right|dt\\
&<&\!\!\frac{|z|^{\mu+\frac{1}{2}}e^{-\Re (z)/2}}{B(\mu\!-\!\kappa\!+\!\fs,\mu\!+\!\kappa\!+\!\fs)}\,\frac{2^\nu|p|^{\nu+1}\g(\nu+\fs)}{\sqrt{\pi} (\Re (p))^{2\nu+1}}
\int_0^1t^{\mu-\kappa+\nu-\frac{1}{2}}(1-t)^{\mu+\kappa+\nu-\frac{1}{2}}e^{\Re (z) t}dt\\
&=&\!\!\frac{2^\nu |p|^{\nu+1} \g(\nu+\fs)}{\sqrt{\pi} (\Re (p))^{2\nu+1}}\,|z|^{\mu+\frac{1}{2}}e^{-\Re (z)/2}\,\frac{B(\mu\!-\!\kappa\!+\!\nu\!+\!\fs,\mu\!+\!\kappa\!+\!\nu\!+\!\fs)}{B(\mu\!-\!\kappa\!+\!\fs,\mu\!+\!\kappa\!+\!\fs)}\\
&&\hspace{6cm}\times \Phi(\mu\!-\!\kappa\!+\!\nu\!+\fs;2\mu\!+\!2\nu\!+\!1;\Re (z)),
\end{eqnarray*}
upon evaluation of the integral by means of (\ref{e3}). This completes the proof.

\begin{theorem}$\!\!\!.$\ Let $p>0$, $\nu\geq0$ and the parameters $\kappa$, $\mu$ be real and such that $\mu\pm\kappa+\fs>0$. Further let $x>0$ and $n$ denote a non-negative integer. Then the following upper bound holds:
\bee\label{e551}
M_{\kappa,\mu}^{(p,\nu)}(x)>x^{\mu+\frac{1}{2}}e^{-x/2} \sum_{j=0}^n \frac{x^j}{j!}\,\frac{B_{p,\nu}(\mu\!-\!\kappa\!+\!j\!+\!\fs,\mu\!+\!\kappa\!+\!\fs)}{B(\mu\!-\!\kappa\!+\!\fs,\mu\!+\!\kappa\!+\!\fs)}.
\ee
\end{theorem}
\noindent{\bf Proof.}\ \ We have for $x>0$ and $t\in(0,1)$
\[e^{xt}>\sum_{j=0}^n \frac{(xt)^j}{j!}\qquad (n=0, 1, 2, \ldots).\]
Then, since $K_{\nu+\frac{1}{2}}(x)>0$ for $x>0$ and $\nu\geq 0$, we have from (\ref{e20})
\[M_{\kappa,\mu}^{(p,\nu)}(x)>\frac{x^{\mu+\frac{1}{2}}e^{-x/2}}{B(\mu\!-\!\kappa\!+\!\fs,\mu\!+\!\kappa\!+\!\fs)}
\sqrt{\frac{2p}{\pi}}\sum_{j=0}^n\frac{x^j}{j!}\int_0^1t^{\mu-\kappa+j-1}(1-t)^{\mu+\kappa-1}K_{\nu+\frac{1}{2}}\bl(\frac{p}{t(1-t)}\br)dt\]
\[=x^{\mu+\frac{1}{2}}e^{-x/2} \sum_{j=0}^n \frac{x^j}{j!}\,\frac{B_{p,\nu}(\mu\!-\!\kappa\!+\!j\!+\!\fs,\mu\!+\!\kappa\!+\!\fs)}{B(\mu\!-\!\kappa\!+\!\fs,\mu\!+\!\kappa\!+\!\fs)}\]
upon use of (\ref{e13}). This completes the proof.

When $n=0$, we obtain the simple lower bound
\bee\label{e552}
M_{\kappa,\mu}^{(p,\nu)}(x)>x^{\mu+\frac{1}{2}}e^{-x/2}\,\frac{B_{p,\nu}(\mu\!-\!\kappa\!+\!\fs,\mu\!+\!\kappa\!+\!\fs)}{B(\mu\!-\!\kappa\!+\!\fs,\mu\!+\!\kappa\!+\!\fs)}.
\ee

\vspace{0.6cm}

\begin{center}
{\bf 6. \  Concluding remarks}
\end{center}
\setcounter{section}{6}
\setcounter{equation}{0}
\renewcommand{\theequation}{\arabic{section}.\arabic{equation}}
We have given the $(p,\nu)$-extension of the Whittaker function of the first kind, which we denote by $M_{\kappa,\mu}^{(p,\nu)}(z)$, by using a similarly extended confluent hypergeometric function introduced in \cite{PCP}. This involves the use of the classical Beta function extended by the introduction in its integral representation of the modified $K$-Bessel function of order $\nu+\frac{1}{2}$ and argument $p/t(1-t)$. Some of the main properties of $M_{\kappa,\mu}^{(p,\nu)}(z)$ have been investigated, namely a differential property, integral representations and the analogue of the Kummer transformation. The leading asymptotic form has been derived for $z\to+\infty$ by application of the saddle-point method. The Mellin transform and an associated integral transform, together with some bounds on $M_{\kappa,\mu}^{(p,\nu)}(z)$, have been obtained. 

The Whittaker function has various applications in mathematical physics and engineering sciences. It is hoped that the results obtained in this paper will have some interest in these fields.


\vspace{0.6cm}

\begin{thebibliography}{11}
\footnotesize{



\bibitem{a2} Andrews, G.E., Askey, R. and Roy, R., \textit{Special Functions}, Cambridge University Press, Cambridge, 1999.

\bibitem{b1} Bailey, W.N., \textit{Generalized Hypergeometric Series}, Cambridge Math. Tract No. 32, Cambridge University Press, Cambridge 1935; Reprinted by Stechert-Hafner, New York, 1964.

\bibitem{B} Chaudhry, M.A., Qadir, A., Rafique, M. and  Zubair, S.M., Extension of Euler's beta function, J. Comput. Appl. Math., {\bf 78} (1997) 19--32.

\bibitem{C} Chaudhry, M.A., Qadir, A., Srivastava, H.M. and Paris, R.B., Extended hypergeometric and confluent hypergeometric functions, Appl. Math. Comput. {\bf 159} (2004) 589--602.

\bibitem{E1} Choi, J., Rathie, A.K. and Parmar, R.K. Extension of extended Beta, hypergeometric function and confluent hypergeometric function, Honam Mathematical J., {\bf 36} (2014) 357-385.

\bibitem{E2} Choi, J., Ghayasuddin, M. and Khan, N., Generalized Extended Whittaker Function and Its Properties, Applied Mathematical Sciences, {\bf 9} (2015) 6529--6541.

\bibitem{DP1} Dar, S.A. and Paris, R.B.,  A $(p,\nu)$-extension of the Appell function $F_1(\cdot)$ and its properties. (2017) arXiv:1711.07780. [submitted for publication].

\bibitem{DS} Dar, S.A. and Shadab, M., $(p,q)$-extension of the Whittaker function and its certain properties. to appear in Korean Math. Soc. (2018).




\bibitem{M} Nagar, D.K., V\'{a}squez, R.A.M. and Gupta, A.K.,  Properties of the
extended Whittaker function, Progr. Appl. Math., {\bf 6}, No. 2 (2013) 70--80.

\bibitem{DLMF}Olver F.W.J., Lozier D.W., Boisvert R.F. and Clark C.W. (eds.), \textit{NIST Handbook of
Mathematical Functions}, Cambridge University Press, Cambridge, 2010.

\bibitem{O1} $\ddot{\mbox{O}}$zergin, E.,  $\ddot{\mbox{O}}$zarslan, M.A. and Altn, A., Extension of gamma, beta and hypergeometric functions, Journal of Computational and Applied Mathematics, {\bf 235} (2011) 4601--4610.

\bibitem{PCP} Parmar, R.K., Chopra, P. and Paris, R.B., On an Extension of Extended Beta and Hypergeometric Functions, J. Classical Anal. {\bf 11} (2017) 91--106. [arXiv:1502.06200].




\bibitem{W} Whittaker, E.T. and Watson, G.N., {\it A Course of Modern Analysis}. Cambridge University Press, Cambridge, 1952.
}
\end{thebibliography}
\end{document}